\documentclass[11pt]{article}
\usepackage{amsmath}
\usepackage{amsmath,amssymb,latexsym,color}
\usepackage[mathscr]{eucal}
\newtheorem{thm}{Theorem}[section]
\newtheorem{defin}[thm]{Definition}

\newtheorem{lemma}[thm]{Lemma}
\newtheorem{cor}[thm]{Corollary}

\newcommand{\proof}{{\it Proof.\quad}}
\newcommand{\qed}{\hfill\Box\medskip}
\usepackage{CJK}
\begin{document}
\begin{CJK*}{GBK}{song}
\renewcommand{\abovewithdelims}[2]{
\genfrac{[}{]}{0pt}{}{#1}{#2}}

\title{\bf On the fractional metric dimension of graphs}

\author{Min Feng\quad Benjian Lv\quad Kaishun Wang\footnote{Corresponding author. E-mail address: wangks@bnu.edu.cn}\\
{\footnotesize   \em  Sch. Math. Sci. {\rm \&} Lab. Math. Com. Sys.,
Beijing Normal University, Beijing, 100875,  China} }
 \date{}
 \maketitle

\begin{abstract}
 In [S. Arumugam, V. Mathew and J. Shen, On fractional metric
dimension of graphs, preprint], Arumugam et al. studied the
fractional metric dimension of the cartesian product of two graphs,
and proposed four open problems. In this paper, we determine the
fractional metric dimension of   vertex-transitive graphs,  in
particular, the fractional metric dimension of a vertex-transitive
distance-regular graph  is expressed in terms of its intersection
numbers. As an application, we calculate the fractional metric
dimension of Hamming graphs and Johnson graphs, respectively.
Moreover, we give an inequality for metric dimension and fractional
metric dimension of an arbitrary graph, and determine all graphs
when the equality holds. Finally, we establish bounds on the
fractional metric dimension of the cartesian product of graphs. As a
result, we completely solve the four open problems.


\end{abstract}

\bigskip

\bigskip

\section{Introduction}

Let $G$ be a finite, simple  and connected graph. We often denote by
$V(G)$ and $E(G)$ the vertex set and the edge set of $G$,
respectively. For any two vertices $x$ and $y$ of $G$, $d_G(x,y)$
denotes the distance between $x$ and $y$, $R_G\{x,y\}$ denotes the
set of vertices $z$ such that $d_G(x,z)\neq d_G(y,z)$. If the graph
$G$ is clear from the context, $d_G(x,y)$ and $R_G\{x,y\}$ will be
written $d(x,y)$ and $R\{x,y\}$, respectively. A {\em resolving set}
of $G$ is a subset $W$ of $V(G)$ such that $W\cap
R_G\{x,y\}\neq\emptyset$ for any two distinct vertices $x$ and $y$
of $G$. The \emph{metric dimension} of $G$, denoted by $\dim(G)$, is
the minimum cardinality of all the resolving sets of $G$. Metric
dimension was first introduced in the 1970s, independently by Harary
and Melter \cite{Ha} and by Slater \cite{Sl}. It is a parameter that
has appeared in various applications (see \cite{RP,Ca} for more
information).

Let $f$: $V(G)\rightarrow[0,1]$ be a real value function. For
$W\subseteq V(G)$, denote $f(W)=\sum_{v\in W}f(v)$. We call $f$ a
{\em resolving function} of $G$ if $f(R_G\{x,y\})\geq 1$ for any two
distinct vertices $x$ and $y$ of $G$. The {\em fractional metric
dimension}, denoted by $\dim_f(G)$, is given by
$$
\dim_f(G)=\min\{|g|:g\textup{ is a resolving function of }G\},
$$
where $|g|=g(V(G))$. Arumugam and Mathew \cite{Ar} formally
introduced the fractional metric dimension of graphs and made some
basic results.

The {\em cartesian product} of graphs $G$ and $H$, denoted by $G\Box
H$, is the graph with the vertex set $V(G)\times V(H)=\{(u, v)|u\in
V(G),v\in V(H)\}$, where $(u_1, v_1)$ is adjacent to $(u_2,v_2)$
whenever $u_1=u_2$ and $\{v_1,v_2\}\in E(H)$, or $v_1=v_2$ and
$\{u_1,u_2\}\in E(G)$. When there is no confusion the vertex $(u,v)$
of $G\Box H$ will be written $uv$. Observe that $d_{G\Box
H}(u_1v_1,u_2v_2)=d_G(u_1,u_2)+d_H(v_1,v_2)$.

Very recently, Arumugam et al. \cite{Aru} characterized all  graphs
$G$ satisfying $\dim_f(G)=\frac{|V(G)|}{2}$, presented several
results on the fractional metric dimension of the cartesian product
 of graphs, and
  raised the following four open problems:

\medskip

{\bf Problem 1.} Determine $\dim_f(K_2\Box C_n)$ when $n$ is odd,
where $K_2$ is the complete graph of order 2 and $C_n$ is a cycle of
order $n$.

{\bf Problem 2.} Determine $\dim_f(H_{n,k})$, where the Hamming graph
$H_{n,k}$ is the cartesian product of $n$ cliques $K_k$.

{\bf Problem 3.} C\'aceres et al. \cite{Ca}  proved $\dim(G\Box
H)\geq\max\{\dim(G),\dim(H)\}.$ Is a similar result true for
$\dim_f(G\Box H)$?

{\bf Problem 4.} Let $G$ and $H$ be two graphs with
$\dim_f(G)=\frac{|V(G)|}{2}$ and $|V(H)|\leq|V(G)|$. Is
$\dim_f(H\Box G)=\frac{|V(G)|}{2}$?

The motivation of this paper is to solve all these problems. In
Section 2, we determine the fractional metric dimension of
vertex-transitive graphs, in particular, the fractional metric
dimension of a vertex-transitive distance-regular graph  is
expressed in terms of its intersection numbers. As an application,
we calculate the fractional metric dimension of Hamming graphs and
Johnson graphs, respectively. In Section 3, we give an inequality
for metric dimension and fractional metric dimension of an arbitrary
graph, and determine all graphs when the equality holds. In Section
4, we establish bounds on the fractional metric dimension of the
cartesian product of graphs.

\section{Vertex-transitive graphs}

For a  graph $G$, in this paper we always assume that
\begin{equation}\label{rG}
r(G)=\min\{|R\{x,y\}|\mid
x,y\in V(G), x\neq y\}.
\end{equation}
In this section we shall express the
fractional metric dimension of a vertex-transitive graph $G$ in
terms of the parameter $r(G)$, and solve Problems 1 and 2.

\begin{lemma}\label{1}
Let $G$ be a graph with $r(G) $ as in {\rm(\ref{rG})}. Then
$\dim_f(G)\leq\frac{|V(G)|}{r(G)}$.
\end{lemma}
\proof Define $f: V(G)\longrightarrow[0,1]$,
$x\longrightarrow\frac{1}{r(G)}$. For any two distinct vertices $x$
and $y$, we have
$$f(R\{x,y\})=\frac{|R\{x,y\}|}{r(G)}\geq 1,$$
which implies that $f$ is a resolving function. Hence,
$\dim_f(G)\leq|f|=\frac{|V(G)|}{r(G)}$.$\qed$

A graph $G$ is {\em vertex-transitive} if its full automorphism
group ${\rm Aut} (G)$ acts transitively on the vertex set.
\begin{thm}\label{vt}
Let $G$ be a vertex-transitive graph with $r(G) $ as in {\rm(\ref{rG})}.
Then $\dim_f(G)=\frac{|V(G)|}{r(G)}$.
\end{thm}
\proof Denote $r=r(G)$.  Then there exist two distinct vertices $u$
and $v$ such that $|R\{u,v\}|=r$. Suppose
$R\{u,v\}=\{w_1,\ldots,w_r\}$.
 For any automorphism $\sigma$ of $G$,
$$R\{\sigma(u),\sigma(v)\}=\{\sigma(w_1),\ldots,\sigma(w_r)\}.$$
Let   $f$ be a resolving function with $\dim_f(G)=|f|$. Then
$$f(\sigma(w_1))+\cdots+f(\sigma(w_r))=f(R\{\sigma(u),\sigma(v)\})\geq1,$$
which implies that
$$\sum_{\sigma\in{\rm Aut} (G)}(f(\sigma(w_1))+\cdots+f(\sigma(w_r)))\geq|{\rm Aut} (G)|.$$
Since $G$ is vertex transitive, we have
$$|{\rm Aut} (G)_{w_1}|\cdot|f|+\cdots+|{\rm Aut} (G)_{w_r}|\cdot|f|\geq|{\rm Aut} (G)|.$$
It follows that $\dim_f(G)=|f|\geq\frac{|V(G)|}{r}$. By Lemma
\ref{1}  we accomplish our proof.$\qed$

Arumugam et al. \cite{Aru} proved that $\dim_f(K_2\Box C_n)=2$ when
$n$ is even. Here we consider the remaining case.

\begin{thm}
If $n$ is an odd integer with $n\geq 3$, then $\dim_f(K_2\Box C_n)=\frac{2n}{n+1}$.
\end{thm}
\proof  For any two  distinct vertices $u_1v_1$ and $u_2v_2$ of
$K_2\Box C_n$, we have
\begin{equation*}
|R\{u_1v_1,u_2v_2\}|=\left\{
\begin{array}{ll}
2n-2,&\textup{ if }u_1=u_2,v_1\neq v_2,\\
2n,&\textup{ if }u_1\neq u_2,v_1=v_2,\\
n+1,&\textup{ if }u_1\neq u_2,d_{C_n}(v_1,v_2)=1,\\
2n-2,&\textup{ if }u_1\neq u_2,d_{C_n}(v_1,v_2)\geq 2.
\end{array}\right.
\end{equation*}
Since $K_2\Box C_n$ is vertex-transitive, $\dim_f(K_2\Box
C_n)=\frac{2n}{n+1}$ by Theorem \ref{vt}. $\qed$

 Next we shall consider the
fractional metric dimension of  distance-regular graphs, in
particular we compute this parameter of Hamming graphs and Johnson
graphs, respectively.

A graph $G$  with diameter $d$  is said to be {\em distance-regular}
if, for all integers $0\leq h,i,j\leq d$ and any two vertices  $x,
y$ at distance $h$, the number
 $$p^h_{i,j} = |\{z\in V(G)\mid d(x,z)=i, d(y,z)=j\}|$$
is a constant. The numbers $ p^h_{i,j}$ are called the {\em
intersection numbers} of $G$. For more information about
distance-regular graphs, we would like to refer readers to
\cite{Br}.

\begin{thm}\label{drg}
Let $G$ be a vertex-transitive distance-regular graph with diameter
$d$. Then
$$\dim_f(G)=\frac{|V(G)|}{|V(G)|-\max\{\sum_{i=1}^dp_{i,i}^h|h=1,\ldots,d\}}.$$
\end{thm}
\proof For any two  distinct vertices $x$ and $y$  at distance $h$,
  $|R\{x,y\}|=|V(G)|-\sum_{i=1}^dp_{i,i}^h$. By Theorem \ref{vt},
the  desired result follows.$\qed$

The {\em Hamming graph}, denoted by $H_{n,k}$, has the vertex set
$\{(x_1,\ldots,x_n)|1\leq x_i\leq k, 1\leq i\leq n\}$, with two
vertices being adjacent if they differ in exactly one co-ordinate.
It is well-known that $H_{n,k}$ is a vertex-transitive
distance-regular graph of order $k^n$, and two vertices are at
distance $j$ if and only if they differ in exactly $j$ co-ordinates.
The {\em hypercube} $Q_n$ is the Hamming graph $H_{n,2}$. Arumugam
and Mathew \cite{Ar} proved  $\dim_f(Q_n)=2$ for $n\geq2$. Now we
compute $\dim_f(H_{n,k})$.

\begin{thm}\label{hamming}
Let $H_{n,k}$   be a Hamming graph where   $k\geq 3$. Then
$\dim_f(H_{n,k})=\frac{k}{2}.$
\end{thm}
\proof Let $\delta_{i,j}$ denote the Kronecker delta.  Pick two
vertices
 $$u=(1,\ldots,1),\quad
 v=(\underbrace{2,\ldots,2}_h,1,\ldots,1).$$
  Then $d(u,v)=h$. Since,
for any vertex $w=(w_1,\ldots,w_n)$, $d(u,w)=d(v,w)$ if and only if
$\sum_{i=1}^h\delta_{1,w_i}=\sum_{i=1}^h\delta_{2,w_i}$, then the
intersection numbers of $H_{n,k}$ satisfy
\begin{eqnarray}\label{ha4}
\sum_{i=1}^{n}p_{i,i}^h=\sum_{s=0}^{\lfloor\frac{h}{2}\rfloor}{h\choose 2s}
{2s\choose s}(k-2)^{h-2s}k^{n-h}.
\end{eqnarray}
Since $\sum_{i=1}^{n}p_{i,i}^1=(k-2)k^{n-1}$, by Theorem \ref{drg}
it suffices to show that
$$
\sum_{i=1}^{n}p_{i,i}^h\leq (k-2)k^{n-1},\quad 2\leq h\leq n.
$$

For $1\leq s\leq\frac{h}{2}$, we have
\begin{eqnarray*}
{h\choose 2s}{2s\choose s}&=&{h-1\choose 2s-1}{2s\choose s}+{h-1\choose 2s}{2s\choose s}\\
&\leq&{h-1\choose 2s-1}\cdot 2^{2s-1}+{h-1\choose 2s}\cdot 2^{2s},
\end{eqnarray*}
which implies  that
\begin{eqnarray*}
&&\sum_{s=0}^{\lfloor\frac{h}{2}\rfloor}{h\choose 2s}{2s\choose
s}(k-2)^{-2s}\\
&\leq&1+\sum_{s=1}^{\lfloor\frac{h}{2}\rfloor}\bigg({h-1\choose 2s-1}\big(\frac{2}{k-2}\big)^{2s-1}+{h-1\choose 2s}\big(\frac{2}{k-2}\big)^{2s}\bigg)\\
&=&\big(1+\frac{2}{k-2}\big)^{h-1}.
\end{eqnarray*}
By (\ref{ha4}), we get
\begin{eqnarray*}
\sum_{i=1}^{n}p_{i,i}^h\leq\big(1+\frac{2}{k-2}\big)^{h-1}(k-2)^hk^{n-h}
=(k-2)k^{n-1},
\end{eqnarray*}
 as
desired. $\qed$

Let $X$ be a set of size $n$, and let ${X\choose k}$ denote the set
of all  $k$-subsets of $X$.
 The {\em Johnson graph}, denoted
by $J(n,k)$, has ${X\choose k}$ as the vertex set, where two
$k$-subsets are adjacent if their intersection has size $k-1$. As we
know, $J(n,k)$ is a vertex-transitive distance-regular graph of
order ${n\choose k}$, and two vertices are at distance $j$ if and
only if their intersection has size $k-j$. Since $J(n,k)\simeq
J(n,n-k)$ and $J(n,1)\simeq K_n$, we only consider the case $4\leq
2k\leq n$. In order to  calculate $\dim_f(J(n,k))$, we need the
following result, the proof of which is immediate from the
unimodality of binomial coefficients.

\begin{lemma}\label{j2}
Let $m$ be a positive integer and $n$ be an arbitrary integer. Then
\begin{eqnarray*}
{m\choose n+1}+{m\choose n-1}\geq{m\choose n}.
\end{eqnarray*}
\end{lemma}

\begin{thm}
Let $J(n,k)$ be a Johnson graph with $4\leq 2k\leq n$. Then
\begin{eqnarray*}
\dim_f(J(n,k))=\left\{
\begin{array}{ll}
3,& \textup{ if }(n,k)=(4,2),\\
\frac{35}{17},& \textup{ if }(n,k)=(8,4),\\
\frac{n^2-n}{2kn-2k^2},& \textup{ otherwise. }
\end{array}\right.
\end{eqnarray*}
\end{thm}
\proof Pick  any two distinct vertices $A$ and $B$ at distance $h$,
write   $A_1=A\setminus (A\cap B)$ and $B_1=B\setminus (A\cap B)$.
Then $|A_1|=|A_2|=h$ and $A_1\cap B_1=\emptyset$. Since, for any
vertex $C$, $d(A,C)=d(B,C)$ if and only if $ |A_1\cap C|=|B_1\cap
C|,$ then the intersection numbers of $J(n,k)$ satisfy
\begin{eqnarray}\label{j6}
\sum_{i=1}^kp_{i,i}^h=\sum_{s=0}^{h}{h\choose s}^2{n-2h\choose
k-2s}.
\end{eqnarray}

If $(n,k)=(4,2)$, by Theorem \ref{drg}  we have $\dim_f(J(4,2))=3$.
If $(n,k)=(8,4)$, by (\ref{j6}) and Theorem \ref{drg} we obtain
$\dim_f(J(8,4))=\frac{35}{17}$.

Now suppose $(n,k)\not\in\{(4,2),(8,4)\}$. Since
$\sum_{i=1}^kp_{i,i}^1={n-2\choose k}+{n-2\choose k-2},$ by
Theorem~\ref{drg} it suffices to show that for  $2\leq h\leq k\leq\frac{n}{2}$,
\begin{eqnarray}
\sum_{i=1}^kp_{i,i}^h \leq{n-2\choose k}+{n-2\choose k-2}.\label{j34}
\end{eqnarray}

We divide our proof into two cases.

{\em Case 1.} $h=\frac{n}{2}$.  Then $h=k$. By (\ref{j6}), we have
\begin{eqnarray*}
\sum_{i=1}^kp_{i,i}^h=\left\{
\begin{array}{ll}
0,&k\textup{ is odd,}\\
{k\choose\frac{k}{2}}^2,&k\textup{ is even.}
\end{array}\right.
\end{eqnarray*}
Since
\begin{eqnarray*}\label{j4}
{k\choose\frac{k}{2}}^2\leq 2{2k-2\choose k-2}={n-2\choose k}+{n-2\choose k-2}\quad
 \textup{ for }k\geq 6,
\end{eqnarray*}
then (\ref{j34}) holds.

{\em Case 2.} $2\leq h<\frac{n}{2}$.   For $1\leq
s\leq h-1$ and $0\leq j\leq 2$, we have
\begin{equation}\label{j5}
{2h-2\choose 2s-j}=\sum_{i=0}^{h-2}{h-2\choose i}{h\choose
2s-j-i}\geq{h-2\choose s-j}{h\choose s}+{h-2\choose s-1}{h\choose
s-j+1}.
\end{equation}
By Lemma \ref{j2} and (\ref{j5}), we have
\begin{eqnarray*}
&&\sum_{i=1}^{2h-2}{2h-2\choose i}{n-2h\choose k-i}+\sum_{i=0}^{2h-3}{2h-2\choose i}{n-2h\choose k-i-2}\\
&=&\sum_{s=1}^{h-1}\bigg\{\bigg[{2h-2\choose 2s}+{2h-2\choose
2s-2}\bigg]{n-2h\choose k-2s}\\
&&\quad+{2h-2\choose 2s-1}\bigg[{n-2h\choose k-2s+1}+{n-2h\choose k-2s-1}\bigg]\bigg\}\\
&\geq&\sum_{s=1}^{h-1}\bigg[{2h-2\choose 2s}+{2h-2\choose 2s-2}+{2h-2\choose 2s-1}\bigg]{n-2h\choose k-2s}\\
&\geq&\sum_{s=1}^{h-1}\bigg[{h-2\choose s}+2{h-2\choose s-1}+{h-2\choose s-2}\bigg]{h\choose s}{n-2h\choose k-2s}\\
&=&\sum_{s=1}^{h-1}{h\choose s}^2{n-2h\choose k-2s}.
\end{eqnarray*}
Then
\begin{eqnarray*}
{n-2\choose k}+{n-2\choose k-2}
 &=&\sum_{i=0}^{2h-2}{2h-2\choose i}{n-2h\choose k-i}+\sum_{i=0}^{2h-2}{2h-2\choose i}{n-2h\choose k-i-2}\nonumber\\
&\geq&{n-2h\choose k}+\sum_{s=1}^{h-1}{h\choose s}^2{n-2h\choose k-2s}+{n-2h\choose k-2h}\nonumber\\
&=&\sum_{s=0}^{h}{h\choose s}^2{n-2h\choose k-2s}.
\end{eqnarray*}
Hence (\ref{j34}) holds by (\ref{j6}).
 $\qed$

\section{An inequality}

In this section, we give an inequality for metric dimension and
fractional metric dimension of any graph, and determine all graphs
when the equality holds.

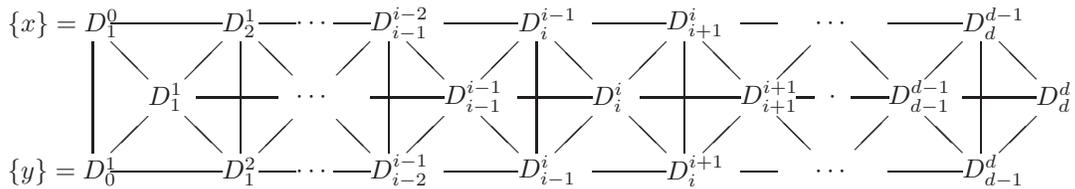
\begin{figure}[htbp]
\unitlength 14pt
\begin{flushright}
\begin{picture}(27,6)
\small
 \put(-1.8,4){\makebox(1,1)[l]{$\{x\}=D^0_1$}}
 \put(4,4){\makebox(1,1)[l]{$D^1_2$}}
 \put(6,4){\makebox(1,1)[l]{$\cdots$}}
 \put(8,4){\makebox(1,1)[l]{$D^{i-2}_{i-1}$}}
 \put(12,4){\makebox(1,1)[l]{$D^{i-1}_{i}$}}
 \put(16,4){\makebox(1,1)[l]{$D^{i}_{i+1}$}}
 \put(19,4){\makebox(3,1){$\cdots$}}
 \put(24,4){\makebox(1,1)[l]{$D^{d-1}_d$}}
 \put(-1.8,0){\makebox(1,1)[l]{$\{y\}=D^1_0$}}
 \put(4,0){\makebox(1,1)[l]{$D^2_1$}}
 \put(6,0){\makebox(1,1)[l]{$\cdots$}}
 \put(8,0){\makebox(1,1)[l]{$D^{i-1}_{i-2}$}}
 \put(12,0){\makebox(1,1)[l]{$D^{i}_{i-1}$}}
 \put(16,0){\makebox(1,1)[l]{$D^{i+1}_{i}$}}
 \put(19,0){\makebox(3,1){$\cdots$}}
 \put(24,0){\makebox(1,1)[l]{$D^{d}_{d-1}$}}
 \put(2,2){\makebox(1,1)[l]{$D^1_1$}}
 \put(6,2){\makebox(1,1)[l]{$\cdots$}}
 \put(10,2){\makebox(1,1)[l]{$D^{i-1}_{i-1}$}}
 \put(14,2){\makebox(1,1)[l]{$D^{i}_{i}$}}
 \put(18,2){\makebox(1,1)[l]{$D^{i+1}_{i+1}$}}
 \put(20,2){\makebox(1,1){$\cdot$}}
 \put(22,2){\makebox(1,1)[l]{$D^{d-1}_{d-1}$}}
 \put(26,2){\makebox(1,1)[l]{$D^d_d$}}
 \put(1,4.5){\line(1,0){3}}
 \put(5,4.5){\line(1,0){1}}
 \put(7,4.5){\line(1,0){1}}
 \put(10,4.5){\line(1,0){2}}
 \put(14,4.5){\line(1,0){2}}
 \put(18,4.5){\line(1,0){1}}
 \put(22,4.5){\line(1,0){2}}
 \put(1,0.5){\line(1,0){3}}
 \put(5,0.5){\line(1,0){1}}
 \put(7,0.5){\line(1,0){1}}
 \put(10,0.5){\line(1,0){2}}
 \put(14,0.5){\line(1,0){2}}
 \put(18,0.5){\line(1,0){1}}
 \put(22,0.5){\line(1,0){2}}
 \put(3.3,2.5){\line(1,0){2.2}}
 \put(8,2.5){\line(1,0){2}}
 \put(12,2.5){\line(1,0){2}}
 \put(15.3,2.5){\line(1,0){2.7}}
 \put(19.5,2.5){\line(1,0){0.5}}
 \put(21,2.5){\line(1,0){1}}
 \put(24,2.5){\line(1,0){2}}
 \put(1,1){\line(1,1){1}}
 \put(13,1){\line(1,1){1}}
 \put(17,1){\line(1,1){1}}
 \put(22,2){\line(-1,-1){0.9}}
 \put(25,1){\line(1,1){1}}
 \put(16,1){\line(-1,1){1}}
 \put(19,2){\line(1,-1){0.9}}
 \put(24,1){\line(-1,1){1}}
 \put(1,4){\line(1,-1){1}}
 \put(13,4){\line(1,-1){1}}
 \put(17,4){\line(1,-1){1}}
 \put(22,3){\line(-1,1){1}}
 \put(25,4){\line(1,-1){1}}
 \put(16,4){\line(-1,-1){1}}
 \put(19,3){\line(1,1){0.9}}
 \put(24,4){\line(-1,-1){1}}
 \put(0.5,1){\line(0,1){3}}
 \put(16.5,1){\line(0,1){3}}
 \put(24.5,1){\line(0,1){3}}
 \put(12.5,1){\line(0,1){3}}
 \put(3,2){\line(1,-1){1}}
 \put(3,3){\line(1,1){1}}
 \put(4.5,1){\line(0,1){3}}
 \put(5,1){\line(1,1){0.9}}
 \put(5,4){\line(1,-1){0.9}}
 \put(8,1){\line(-1,1){0.9}}
 \put(8,4){\line(-1,-1){0.9}}
 \put(9,1){\line(1,1){1}}
 \put(9,4){\line(1,-1){1}}
 \put(8.5,1){\line(0,1){3}}
 \put(12,1){\line(-1,1){1}}
 \put(12,4){\line(-1,-1){1}}
\end{picture}
\end{flushright}
\caption{\label{diagram:rank1} Intersection diagram with respect to
$x$ and $y$.}
\end{figure}

\begin{lemma}\label{ex0}
Let $G$ be a graph. For any   subset $A$ of $V(G)$ with size
$|V(G)|-\dim(G)+1$,   there exist two distinct vertices $x$ and $y$
of $G$ such that $R\{x,y\}\subseteq A$.
\end{lemma}
\proof Suppose there exists a subset $A$ with size
$|V(G)|-\dim(G)+1$ such that $R\{x,y\}\not\subseteq A$ for any two
distinct vertices $x$ and $y$. Then $R\{x,y\}\cap (V(G)\setminus
A)\neq\emptyset$; and so  $V(G)\setminus A$ is a resolving set of
$G$. Therefore, $\dim(G)-1=|V(G)\setminus A|\geq\dim(G)$,   a
contradiction. $\qed$
\begin{lemma}\label{ex1}
Let $G$ be a graph with $r(G) $ as in {\rm(\ref{rG})}. Then
$r(G)=|V(G)|-1$ if and only if $G$ is isomorphic to a path or an odd
cycle.
\end{lemma}
\proof The sufficiency is immediate. Conversely, suppose
$r(G)=|V(G)|-1$. Denote the maximum degree of $G$ by $\Delta$. Pick
a vertex $x$ with degree  $\Delta$. Suppose $\Delta\geq3$. We may
choose three pairwise distinct vertices $x_1, x_2$ and $x_3$
adjacent to $x$. Observe $x\not\in R\{x_1,x_2\}\cup R\{x_1,x_3\}\cup
R\{x_2,x_3\}$. Then $x_3\in R\{x_1,x_2\}$, which implies that
$d(x_3,x_1)\neq d(x_3,x_2)$. We may assume $d(x_3,x_1)=1$ and
$d(x_3,x_2)=2$. From $x_2\in R\{x_1,x_3\}$ we get $d(x_1,x_2)=1$,
which implies that $x_1\not\in R\{x_2,x_3\}$, a contradiction. Hence
$\Delta\leq2$; and so $G$ is isomorphic to a path or a cycle. If $G$
is isomorphic to an even cycle, then for two vertices $u$ and $v$ at
distance 2, we have $|R\{u,v\}|=n-2$, a contradiction. Hence, the
desired result follows. $\qed$
\begin{lemma}\label{ex2}
Let $G$ be a graph. Suppose $|\overline{R}\{u,v\}|=2$ for any two
distinct vertices $u$ and $v$, where  $\overline
R\{u,v\}=V(G)\setminus{R}\{u,v\}$.  If the map
$$\varphi: {V(G)\choose2}\longrightarrow{V(G)\choose2}, \quad\{u,v\}\longmapsto \overline
R\{u,v\}$$ is a bijection, then $G$ is isomorphic to the complete
graph of order four.
\end{lemma}
\proof For any two adjacent  vertices $x$ and $y$, let
$D^i_j(x,y)=\{u\in V(G)\mid d(x,u)=i,d(y,u)=j\}$. The intersection
diagram with respect to $x$ and $y$ is the collection
$\{D^i_j(x,y)\}_{i,j}$ with lines between $D^i_j(x,y)$'s and
$D^s_t(x,y)$'s. We draw a line between $D^i_j(x,y)$ and $D^s_t(x,y)$
if there is possibility of existence of edges. The intersection
diagram with respect to $x$ and $y$ is shown in Figure
\ref{diagram:rank1}, where $D^i_j=D^i_j(x,y)$ and $d$ is the
diameter of $G$.

Since $\overline{R}\{x,y\}=\cup_{i=1}^dD^i_i$, then
\begin{equation}\label{exi2}
\sum_{i=1}^d|D^i_i|=2.
\end{equation}

Note that two distinct vertices $u$ and $v$ belong to some $D^i_j$
if and only if $\overline{R}\{u,v\}=\{x,y\}$. Since $\varphi$ is a
bijection, there exist $a$ and $b$ such that $|D^a_b|=2$ and
\begin{eqnarray}
|D^i_j|\leq 1\mbox{ for }(i,j)\neq (a,b).\label{exi4}
\end{eqnarray}
Write  $D^a_b=\{z_1,z_2\}$.

{\bf Claim.} There exist two adjacent vertices $x_0$ and $y_0$
such that $|D^1_1(x_0,y_0)|=2$.

  Suppose that, for any two adjacent vertices $x$ and $y$,
\begin{equation}\label{exi5}
|D^1_1(x,y)|\leq 1.
\end{equation}

{\em Case 1.} $a=b$. By (\ref{exi2}), $D^i_i=\emptyset$ for $i\neq a$. Then
there exist
$x'\in D^{a-1}_{a}$ and
 $y'\in D^{a}_{a-1}$ such that
$\overline{R}\{x',y'\}=\{z_1,z_2\}=\overline{R}\{x,y\}$. Hence
$\{x',y'\}=\{x,y\}$ and $|D^1_1|=2$,   contrary to (\ref{exi5}).

{\em Case 2.} $a\neq b$. By (\ref{exi4}),
there exists a unique  $z\in D^{a-1}_{b-1}$ such that
$z\in \overline{R}\{z_1,z_2\}=\{x,y\}$. Therefore, $(a,b)=(1,2)$ or $(2,1)$.
 We may assume $D^1_2=\{z_1,z_2\}$.
 By (\ref{exi2}) and (\ref{exi4}), there exist
$i_0<j_0$ such that $|D^{i_0}_{i_0}|=|D^{j_0}_{j_0}|=1$ and
$D^i_i=\emptyset$ for $i\neq i_0,j_0$. Write $D^{i_0}_{i_0}=\{w_1\}$
and $D^{j_0}_{j_0}=\{w_2\}$. Since
$w_1\not\in\overline{R}\{z_1,z_2\}$, then $ d(w_1,z_1)\neq
d(w_1,z_2). $ Without loss of generality,  assume that
$d(w_1,z_1)<d(w_1,z_2).$

{\em Case 2.1.} $i_0=1$. Then $d(w_1,z_1)=1$, which implies
$D^1_1(x,w_1)=\{y,z_1\}$. Consequently, $x, w_1$ are adjacent and
$|D^1_1(x,w_1)|=2$,  contrary to (\ref{exi5}).

{\em Case 2.2.} $i_0=2$. Note that $D^1_1=\emptyset$. Then
$d(w_1,z_1)=d(w_1,y_1)=1$, where $y_1$ is the unique vertex in
$D^2_1$. If $D^2_3=\emptyset$ or $D^3_2=\emptyset$, then $j_0=3$ and
$d(w_1,w_2)=1$. Consequently,
$\overline{R}\{z_1,y_1\}=\{w_1,w_2\}=\overline{R}\{x,y\}$, which
contradicts the fact that $\varphi$ is a bijection. Write
$D^2_3=\{x_2\}$ and $D^3_2=\{y_2\}$. Since
$\overline{R}\{z_1,y_1\}\neq\{w_1,w_2\}$, then $d(w_2,z_1)\neq
d(w_2,y_1)$. The fact that $d(w_2,x_2)=d(w_2,y_2)$ implies that
$x_2$ is not adjacent to $z_1$ and $d(x_2,z_2)=1$; and so
$V(G)\setminus\{z_1,w_1,w_2\}\subseteq R\{x_2,y_2\}\cap
R\{z_2,y_1\}$. Then $\overline{R}\{x_2,y_2\} \cup
\overline{R}\{z_2,y_1\}\subseteq\{z_1,w_1,w_2\}$. Since
$w_2\in\overline R\{x_2,y_2\}$ and
$\overline{R}\{x,y\}=\{w_1,w_2\}$, we get
$\overline{R}\{x_2,y_2\}=\{z_1,w_2\}$. Consequently,
$\overline{R}\{z_2,y_1\}=\{z_1,w_1\}$, which implies
$d(w_1,z_2)=d(w_1,y_1)=d(w_1,z_1)$. Hence $w_1\in
\overline{R}\{z_1,z_2\}$, which contradicts
$\overline{R}\{z_1,z_2\}=\{x,y\}$.

{\em Case 2.3.} $i_0\geq 3$. Note that
$D^{i_0-1}_{i_0-1}=\emptyset$. By (\ref{exi4}),
$|D^{i_0-1}_{i_0}|=|D^{i_0}_{i_0-1}|=1$. Write
$D^{i_0-1}_{i_0}=\{x'\}$ and $D^{i_0}_{i_0-1}=\{y'\}$.
Since $d(x',w_1)=d(y',w_1)$ and $d(x',w_2)=d(y',w_2)$,
then $\overline{R}\{x',y'\}=\{w_1,w_2\}=\overline{R}\{x,y\}$,  a
contradiction.

Therefore, our claim is valid.

Now write $D^1_1(x_0,y_0)=\{z_1',z_2'\}$. By (\ref{exi2}),
$D^i_i(x_0,y_0)=\emptyset$ for $i\geq 2$. By (\ref{exi4}),
$|D^i_j(x_0,y_0)|\leq 1$ for $i\neq j.$

Suppose $D^1_2(x_0,y_0)\cup D^2_1(x_0,y_0)\neq\emptyset$. We may
assume that $D^1_2(x_0,y_0)\neq \emptyset$. Write
$D^1_2(x_0,y_0)=\{x_1\}$. If $D^2_1(x_0,y_0)=\emptyset$, then
$V(G)\setminus\{z_1',z_2'\}\subseteq R\{x,x_1\}$, which implies that
$\overline{R}\{x,x_1\}=\{z_1',z_2'\}=\overline{R}\{x,y\}$, a
contradiction. If $D^2_1(x_0,y_0)\neq\emptyset$, write
$D^2_1(x_0,y_0)=\{y_1\}$, then $V(G)\setminus\{z_1',z_2'\}\subseteq
R\{x_1,y_1\}$. Hence, $\overline{R}\{x_1,y_1\}=\{z_1',z_2'\}$, a
contradiction. Consequently, $D^1_2(x_0,y_0)=
D^2_1(x_0,y_0)=\emptyset$, and $|V(G)|=4$. Since
$\overline{R}\{x_0,z_1'\}=\{y_0,z_2'\}$,  we have $d(z_1',z_2')=1$;
and then $G\simeq K_4$. $\qed$

Now we state our main result of this section.
\begin{thm}\label{prop}
Let $G$ be a graph of order $n$. Then
\begin{equation}\label{exi1}
\dim_f(G)\geq\frac{n}{n-\dim(G)+1}.
\end{equation}
Moreover, the equality   holds if and only if $G$ is isomorphic to a
path, a complete graph, or an odd cycle.
\end{thm}
\proof Write $l=n-\dim(G)+1$.   Suppose $f$ is a resolving function
of $G$ with $|f|=\dim_f(G)$. By Lemma \ref{ex0}, $f(A)\geq 1$ for
each $A\in{V(G)\choose l}.$    Hence
$\sum_{A\in{V(G)\choose{l}}}f(A)\geq{n\choose{l}}$.  Since
$\sum_{A\in{V(G)\choose{l}}}f(A)={n-1\choose l-1}|f|$, then
(\ref{exi1}) holds.

Suppose that the equality in (\ref{exi1}) holds. Then $f(A)=1$ for
each $A\in{V(G)\choose l}$. If $\dim(G)=1$, then $G\simeq P_n$. If
$\dim(G)=n-1$, then $G\simeq K_n$. Now suppose $2\leq\dim(G)\leq
n-2$. Then $3\leq l\leq n-1$.

Given two distinct vertices $x,y$, pick an $(l-1)$-subset $A_1$  of
$V(G)\setminus\{x,y\}$. Since $f(\{x\}\cup A_1)=1=f(\{y\}\cup A_1)$,
then $f(x)=f(y)=\frac{1}{l}$, which implies that $|R\{x,y\}|\geq l$.
By Lemma \ref{ex0}, for any $A\in{V(G)\choose l}$, there exist two
distinct vertices $x_0$ and $y_0$ such that $R\{x_0,y_0\}=A$. Hence
$r(G)=l$, and
\begin{eqnarray}
{n\choose l}\leq|\{R\{u,v\}|u\neq v\}|\leq{n\choose 2}\label{ex4}.
\end{eqnarray}
It follows that $l=n-1$ or $l=n-2$.

{\em Case 1.} $l=n-1$. By   Lemma \ref{ex1}, $G$ is isomorphic to
 an odd cycle.

{\em Case 2.} $l=n-2$. In this case $n\geq 5.$ By (\ref{ex4}), we
have $ |\{R\{u,v\}|u\neq v\}|={n\choose l}$. By Lemma \ref{ex0} we
get $|R\{x,y\}|=l$ for any two distinct vertices $x$ and $y$. Then
we obtain a bijection $\varphi$ as in Lemma~\ref{ex2}. Hence,
$G\simeq K_4$, a contradiction.

The converse is true by \cite[Corollary 2.7 and Theorem 3.2]{Ar}.
$\qed$

Combining Lemma \ref{1} and Theorem \ref{prop}, we obtain the
following corollary.
\begin{cor}
Let $G$ be a  distance-regular graph with diameter $d$. Then
$$\dim(G)\leq \max\{\sum_{i=1}^{d}p_{i,i}^h|h=1,\ldots,d\}+1.$$
The equality holds if and only if $G$ is a complete graph or an odd
cycle.
\end{cor}

\section{Cartesian product of graphs}

 In this section,
we shall establish bounds on the fractional metric dimension of the
cartesian product of two graphs.

\begin{thm}\label{low}
Let $G$ and $H$ be two graphs. Then
$\dim_f(G\Box H)\geq\dim_f(G).$
\end{thm}
\proof
Pick a resolving function $f_{G\Box H}$
 of $G\Box H$ with $|f_{G\Box H}|=\dim_f(G\Box
H)$. Define
$$f_G:  V(G)\longrightarrow [0,1],\quad  u \longmapsto \min\{1,\sum_{y\in
V(H)}f_{G\Box H}(uv)\}.$$ Let $u_1$ and $u_2$ be two distinct
vertices of $G$. We shall prove
\begin{equation}\label{l1}
f_{G}(R_G\{u_1,u_2\})\geq 1.
\end{equation}

If there exists $u_0\in R_G\{u_1,u_2\}$ with $f_G(u_0)=1$, then
(\ref{l1}) holds. Now we suppose $f_G(u)=\sum_{v\in V(H)}f_{G\Box
H}(uv)$ for any  $u\in V(G)$. For  $v_0\in V(H)$,
 we have
\begin{equation*}
R\{u_1v_0,u_2v_0\}=\bigcup_{u\in
R_G\{u_1,u_2\}}\bigcup_{v\in V(H)}\{uv\}.
\end{equation*}
Then
$$f_{G}(R_G\{u_1,u_2\})=\sum_{u\in R_G\{u_1,u_2\}}\sum_{v\in V(H)}f_{G\Box
H}(uv)=f_{G\Box H}(R\{u_1v_0,u_2v_0\})\geq 1,$$
(\ref{l1}) holds.
Therefore, $f_G$
is a resolving function of $G$. Since
$$|f_G|\leq\sum_{u\in
V(G)}\sum_{v\in V(H)}f_{G\Box H}(uv)=|f_{G\Box H}|,$$
then
$\dim_f(G)\leq\dim_f(G\Box H)$, as desired.$\qed$

Since $G\Box H$ is isomorphic to $H\Box G$, this theorem gives an
answer to Problem 3. By Theorem \ref{hamming} the bound in Theorem
\ref{low} is sharp.

\begin{thm}\label{upp}
Let $G$ and $H$ be two graphs.
Then
$$\dim_f(G\Box H)\leq\max\{\dim_f(G),|V(H)|\}.$$
\end{thm}
\proof Let $f_G$ be a resolving function of $G$ with $|f_G|=\dim_f(G)$.
Denote $l=\min\{\dim_f(G),|V(H)|\}$.
Define
$$f_{G\Box H}:V(G\Box H)\longrightarrow[0,1],\quad uv\longmapsto\frac{f_G(u)}{l}.$$
For any two distinct vertices $u_1v_1$ and $u_2v_2$ in $G\Box H$, we
shall prove
\begin{equation*}\label{in8}
f_{G\Box H}(R\{u_1v_1,u_2v_2\})\geq 1.
\end{equation*}

{\em Case 1.} $v_1=v_2$. Since
$$R\{u_1v_1,u_2v_1\}=\bigcup_{u\in R_G\{u_1,u_2\}}\bigcup_{v\in V(H)}\{uv\},$$
 then
\begin{eqnarray*}
f_{G\Box H}(R\{u_1v_1,u_2v_1\}) =\sum_{u\in
R_G\{u_1,u_2\}}\sum_{v\in V(H)}\frac{f_G(u)}{l}
=\frac{|V(H)|}{l}\cdot f_G\{R_G\{u_1,u_2\}\} \geq 1.
\end{eqnarray*}

{\em Case 2.} $v_1\neq v_2$. Write
\begin{equation*}
\begin{array}{l}
S_1=\{ u\in V(G)\mid d_G(u_1,u)=d_G(u_2,u)\},\\
S_2=\{ u\in V(G)\mid d_G(u_1,u)<d_G(u_2,u)\},\\
S_3=\{ u\in V(G)\mid d_G(u_1,u)>d_G(u_2,u)\}.
\end{array}
\end{equation*}
Then
$$R\{u_1v_1,u_2v_2\}\supseteq\big(\bigcup_{u\in S_1\cup S_2}\{uv_1\}\big)\cup\big(\bigcup_{u\in S_1\cup S_3}\{uv_2\}\big).$$
It follows that
\begin{eqnarray*}
f_{G\Box H}(R\{u_1v_1,u_2v_2\})&\geq&\sum_{u\in S_1\cup S_2}f_{G\Box H}(uv_1)+\sum_{u\in S_1\cup S_3}f_{G\Box H}(uv_2)\\
&=&\frac{f_G(S_1\cup S_2)}{l}+\frac{f_G(S_1\cup S_3)}{l}\\
&=&\frac{\dim_f(G)}{l}+\frac{f_G(S_1)}{l}\\
&\geq&1.
\end{eqnarray*}

Therefore, $f_{G\Box H}$ is a resolving function of $G\Box H$. Since
\begin{eqnarray*}
|f_{G\Box H}| =\sum_{v\in V(H)}\sum_{u\in V(G)}\frac{f_G(u)}{l}
=\sum_{v\in V(H)}\frac{\dim_f(G)}{l} =\max\{\dim_f(G),|V(H)|\},
\end{eqnarray*}
the desired result follows.$\qed$

By \cite[Theorem 4.2]{Ar} and \cite[Theorem 3.3]{Aru}, $P_n\Box K_2$
and  $C_{2n}\Box K_2$ meet the bound in Theorem~\ref{upp}.

Finally, we focus on Problem 4.

\begin{thm}\label{vg2}
Let $G$ be a graph with at least three vertices and $\dim_f(G)=\frac{|V(G)|}{2}$.
Let $H$ be a graph with $|V(H)|\leq|V(G)|$. Then $\dim_f(G\Box H)=\frac{|V(G)|}{2}$.
\end{thm}
\proof By Theorem \ref{low}, $\dim_f(G\Box H)\geq\frac{|V(G)|}{2}$.
In order  to prove $\dim_f(G\Box H)\leq\frac{|V(G)|}{2}$, by Lemma
\ref{1}  it suffices to show that
\begin{equation}\label{in4}
|R\{u_1v_1,u_2v_2\}|\geq 2|V(H)|
\end{equation}
holds  for any two distinct vertices $u_1v_1$ and $u_2v_2$  in
$G\Box H$.

{\em Case 1.} $u_1=u_2$. Since
$R\{u_1v_1,u_1v_2\}\supseteq\{uv_1|u\in V(G)\}\cup\{uv_2|u\in V(G)\},$
then
$|R\{u_1v_1,u_1v_2\}|\geq 2|V(G)|\geq2|V(H)|,$
 (\ref{in4}) holds.

{\em Case 2.} $u_1\neq u_2$. For $v\in V(H)$, let
$$
S_v=\{u|u\in V(G),d_G(u_1,u)-d_G(u_2,u)\neq k_v\},
$$
where $k_v=d_H(v_2,v)-d_H(v_1,v)$. Note that
$R\{u_1v_1,u_2v_2\}=\cup_{v\in V(H)}\{uv|u\in S_v\}$. In order to
prove (\ref{in4}),  we only need to show that
\begin{equation}\label{in5}
|S_v|\geq2.
\end{equation}

{\em Case 2.1.} $k_v\neq d_G(u_1,u_2)$ and $k_v\neq -d_G(u_1,u_2)$.
Then $u_1,u_2\in S_v$, and (\ref{in5}) holds.

{\em Case 2.2.} $k_v=d_G(u_1,u_2)$. Then $u_1\in S_v$. Since
 $\dim_f(G)=\frac{|V(G)|}{2}$, by \cite[Theorem 2.2]{Aru}
there exists a vertex $u_2'\in V(G)\setminus\{u_2\}$ such that, for
any $u\in V(G)\setminus\{u_2,u_2'\}$,
\begin{equation}\label{in6}
d_G(u_2',u)=d_G(u_2,u).
\end{equation}
If $u_2'\neq u_1$,  by (\ref{in6}) we have
$d_G(u_1,u_2')-d_G(u_2,u_2')<d_G(u_1,u_2')=k_v$, which implies
$u_2'\in S_v$ and (\ref{in5}) holds. Now suppose $u_2'=u_1$.
 Choose $u_3\in V(G)\setminus\{u_1,u_2\}$. By (\ref{in6}),
$d_G(u_1,u_3)-d_G(u_2,u_3)=0<k_v$. Then $u_3\in S_v$, and so
(\ref{in5}) holds.

{\em Case 2.3.} $k_v=-d_G(u_1,u_2)$. Similar to Case 2.2,
(\ref{in5}) holds.$\qed$

\section*{Acknowledgement}
The authors would like to thank Professor S. Arumugam  for sending
them the preprint \cite{Aru}.  This research was supported by NSF of
China, NCET-08-0052, and the Fundamental Research Funds for the
Central Universities of China.

\end{CJK*}

\end{document}